\documentclass[12pt,a4paper]{article}

\usepackage{a4wide,amssymb,amsmath,amsthm,xspace,epsfig}

\begin{document}

\newcommand{\N}{\Bbb N}
\newcommand{\R}{\Bbb R}
\newcommand{\Z}{\Bbb Z}
\newcommand{\Q}{\Bbb Q}
\newcommand{\C}{\Bbb C}
\newcommand{\PP}{\mathbb{P}}

\newcommand{\LL}{\Bbb L}

\newcommand{\esp}{\vskip .3cm \noindent}
\mathchardef\flat="115B

\newcommand{\lev}{\text{\rm Lev}}

\def\ut#1{$\underline{\text{#1}}$}
\def\CC#1{${\cal C}^{#1}$}
\def\h#1{\hat #1}
\def\t#1{\tilde #1}
\def\wt#1{\widetilde{#1}}
\def\wh#1{\widehat{#1}}
\def\wb#1{\overline{#1}}

\def\restrict#1{\bigr|_{#1}}

\newtheorem{thm}{Theorem}
\newtheorem{lemma}[thm]{Lemma}

\newtheorem{defi}[thm]{Definition}
\newtheorem{conj}[thm]{Conjecture}
\newtheorem{cor}[thm]{Corollary}
\newtheorem{prop}[thm]{Proposition}
\newtheorem{prob}{Problem}
\newtheorem*{rem}{Remark}
\newtheorem*{aside}{Aside}

\title{Some steps on short bridges:\\ Non-metrizable surfaces and CW-complexes}
\date{\today}
\author{Mathieu Baillif}

\maketitle

\abstract{Among the classical variants of the Pr\"ufer surface, some are homotopy equivalent to a CW-complex
(namely, a point or a wedge of a continuum of circles) and some are not.
The obstruction comes from the existence of uncountably many `infinitesimal bridges' linking two metrizable open subsurfaces
inside the surface. We show that any non-metrizable surface that possesses such a system of 
infinitesimal bridges cannot
be homotopy equivalent to a CW-complex. 
More than for the result on its own, we were motivated by trying to blend elementary techniques of algebraic and set-theoretic topology.}

\section{Introduction}
The aim of this note is twofolds. First, 
we provide a generalisation of a result of Gabard \cite{GabardWouuuh} who
showed that the double of the bordered Pr\"ufer surface is not homotopy equivalent to a CW-complex.
Second, we have tried (just for the fun of it\footnote{Though we would understand if the reader has another definition of `fun'.}) 
to blend elementary techniques of two fields that do not intersect often, that is,
set-theoretic and algebraic topology. 

It is well known (see Milnor's famous article \cite{MilnorCW}) that all metrizable manifolds are 
homotopy equivalent to a countable CW-complex. 
Here we investigate what can happen when the manifolds become non-metrizable.  
To avoid confusion, let us fix the terminology:
a {\em manifold} is a Hausdorff space locally homeomorphic to some $\R^n$, which we moreover 
assume connected for simplicity.
A {\em manifold with boundary}, or {\em bordered manifold}, is defined similarly, as usual.
It happens that some non-metrizable manifolds are contractible, 
for instance the original collared Pr\"ufer surface (see below), but some
are not homotopy equivalent to a CW-complex.
Non-metrizable manifolds can be roughly divided in two classes: the ones that are `big'
(for instance, the long ray $\LL_+$), whose
non-metrizability comes from the size of the manifold, and the
ones that are `weird' in the sense that they do have a metrizable subspace whose closure is non-metrizable.
This note is about finding a criterion for manifolds in the latter class which impedes them
to have the homotopy type of a CW-complex.
The former class (the so called Type I manifolds) is treated (in part) in 
\cite{meszigues+Nyikos}.

The theorem we shall prove is the following:

\begin{thm}
  \label{thmMSIB}
  A surface containing a full MSIB is not homotopy equivalent to a CW-complex.
\end{thm} 

The definitions of MSIB and the fullness condition are given below, but roughly,
it means that the surface contains two non-metrizable subsurfaces with common boundary, both   
of them having a metrizable interior. It is thus their common boundary that 
makes them non-metrizable, and they can be thought as being linked by a system of 
(uncountably many) infinitesimal bridges (the boundary components), see below.
In Section~\ref{sec2} 
we give a bunch of motivating examples showing that if these bridges are
`long' (i.e., not infinitesimal), then the surface may well be homotopy
equivalent to a CW-complex.

In the course of the proof of Theorem \ref{thmMSIB}, we shall prove the (easy) following proposition, which is probably folklore
but we have not been able to track down
a reference.

\begin{prop}\label{nocircles}
  A connected surface with boundary and metrizable interior has at most countably many
  boundary components that are circles.
\end{prop}

  It is easy to build examples of bordered surfaces with continuum many circle boundaries 
  (for instance, put a circle boundary in each of the   collars of 
  $\PP_C$ defined below). 
  Though,
  a sequentially compact manifold has only finitely many boundary components, for if there are infinitely many, taking a point in each 
  yields a sequence whose accumulation point(s) cannot have a neighborhood homeomorphic to either $\R^n$ or $\R^{n-1}\times\R_{\ge 0}$.

We 
will not need much sophistication in set theoretic topology, a basic knowledge of ordinals and
non-metrizable manifolds should suffice, though a familiarity with the Pr\"ufer 
surface\footnote{which is by the way more a geometric than a set-theoretic object...} would probably
be helpful. 
The algebraic topology techniques we will use are all elementary, in fact they more or less cover
the tools that one encounters in an introductory course on homotopy theory:
van Kampen's Theorem, Whitehead's Theorem, CW-approximation, universal coverings. An excellent reference is
Hatcher \cite{Hatcher}.
We shall need some results that are specific to surfaces, for instance the fact that
a nullhomotopic circle in a surface bounds a $2$-disk, or that $\R^2$ is the only simply connected non-compact metric surface.
We shall sometimes repeat parts of proofs (some of which not due to the author) published elsewhere for completeness.
\\
We shall be quite loose in our use of pointed spaces and homotopy, and never state out explicitely the base points, though they are
always implicit whenever needed.

\section{Definitions and motivation}\label{sec2}

We first treat a bunch of examples (all are already well known, but we shall
repeat the proofs as they do shed lights on the techniques we shall use later on). 
The Pr\"ufer surface is usually defined as
a half plane $\R\times\R_{\ge 0}$ where each point $(x,0)$ on the boundary is first
replaced by an interval, using a system of wedges pointing at $(x,0)$ union subintervals of the (added) interval
as neighborhoods, and then collaring, i.e., gluing half planes 
$\R\times\R_{\le 0}$ on the new boundary intervals created (see
for instance \cite{GabardWouuuh} for a recent complete description). Here, we shall adopt the equivalent definition
found in Calabi and Rosenlicht's paper \cite{CalabiRosenlicht}:
$\PP$ is a union of a continuum of Euclidean planes $E_r$ ($r\in\R$) with coordinates $(x_r,y_r)$ 
quotiented by the equivalence relation
$(x_r,y_r)\sim(x_s,y_s)$ whenever
\begin{align}
  \label{CR1} 
    y_r=y_s,&\quad\text{ and} \\    
  \label{CR2}
    \begin{array}{lr} x_r y_r+r=x_sy_s+s&\text{ if }y_r=y_s>0, \\
    r=s,\quad x_r=x_s&\text{ if }y_r=y_s\le 0.\end{array}
\end{align} 
Using the map $(x_r,y_r)\mapsto x_ry_r+r$ for $y_r> 0$, we obtain the other construction alluded to above.
We then define the following surfaces:
\begin{align*}
   \PP_s & = \{(x_r,y_r)\in\PP\,:\, y_r\ge 0\} & \text{(bordered Pr\"ufer Surface)},\\
   \PP_C & = \{(x_r,y_r)\in\PP\,:\, y_r\ge -1\}& \text{(bordered collared Pr\"ufer Surface)},\\
   \mathbb{M} & = \PP_s/(x_r,0)\sim(-x_r,0)& \text{(Moore Surface)}.\\
\end{align*}
$\PP_s$ is a separable bordered surface,
its boundary components being $I_r=\{(x_r,0)\in\PP\}$ for $r\in\R$, while $\PP_C$ is a `collared' version of $\PP_s$, which is not 
separable because of the presence of these collars. (Recall that a space is separable if it possesses a 
countable dense subset.)
Beware that some authors call $\PP_s$ the Pr\"ufer surface, and not $\PP$.
$\mathbb{M}$ was first described by R.L. Moore, and is a boundaryless separable surface.
The figure below shows the usual way of picturing $\PP$ (the entire surface drawn), $\PP_C$ (where the parts in the darkest grey are dropped except
their upper boundary), and $\PP_s$ (the lightest grey part, with the bottom boundaries).

\begin{center}
\epsfig{figure=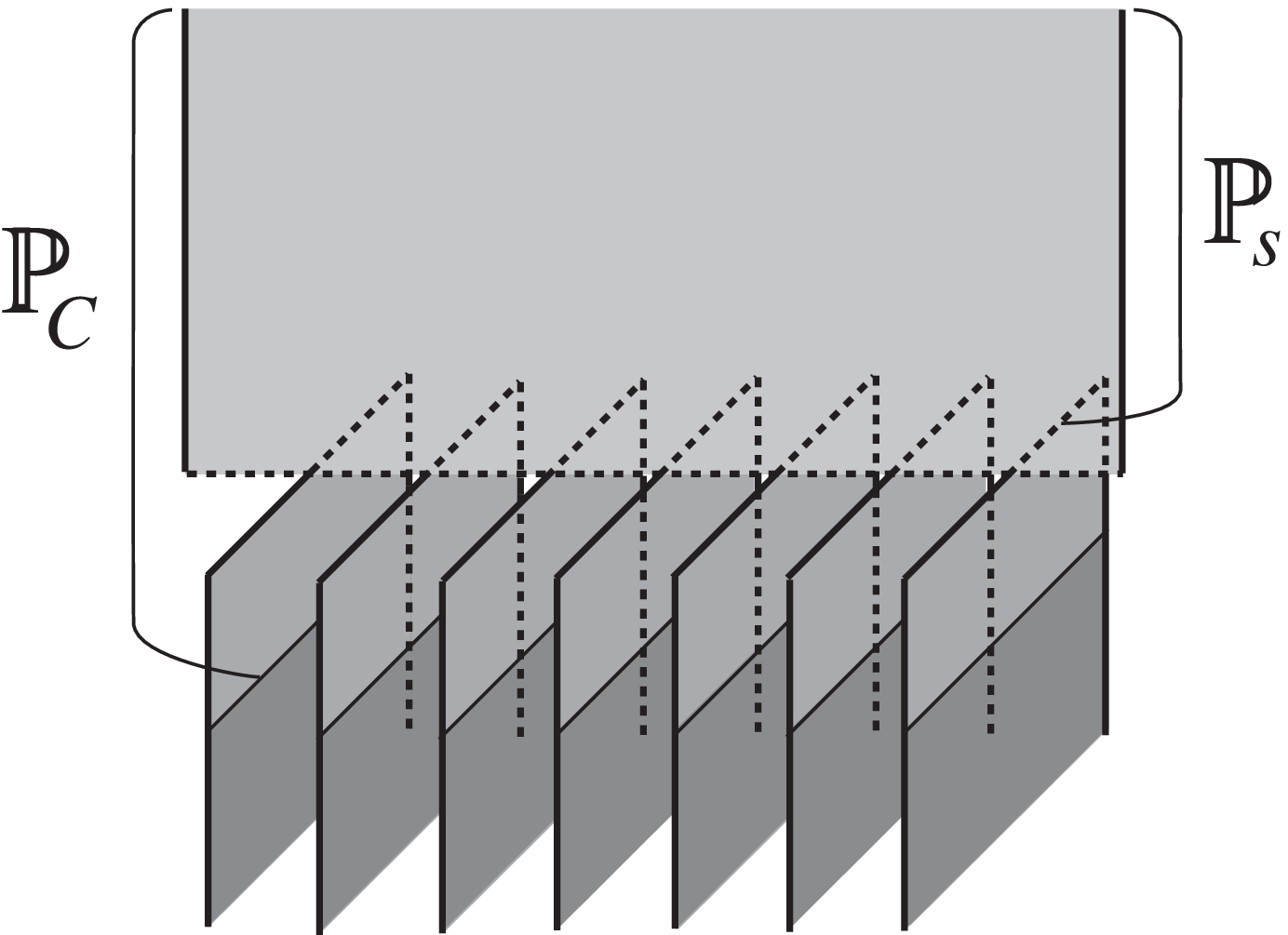, height=4cm}
\end{center}

If $M$ is a bordered manifold, we shall write $2M$ for the double of $M$, i.e. the manifold
without boundary given by two copies $M^0,M^1$ of $M$ whose boundaries are identified pointwise.
The proposition that motivates the results in this note is the following:

\begin{prop}[\cite{CalabiRosenlicht, GabardWouuuh}] \label{propop}
   $\PP$, $\PP_s$, $\PP_C$, $\mathbb{M}$ and $2\PP_C$ are homotopy equivalent to a CW-complex:
   the first four to a point and $2\PP_C$
   to a wedge of continuum many circles;
   whereas $2\PP_s$ is not homotopy equivalent to a CW-complex. 
\end{prop}

\proof[Idea of the proof]  
Consider the following homotopy in $\PP$, taken from
\cite[p. 339]{CalabiRosenlicht}:
\begin{equation}\label{htCR}
h_t(x_r,y_r)=\left\{
  \begin{array}{lr}
     \left(x_r\left(\frac{1-t+ty_r}{1+ty_r}\right)^{1/2}\, , \,
     y_r\left(\frac{1+ty_r}{1-t+ty_r}\right)^{1/2}\right)
     & \text{ when }y_r>0,  \\
     (x_r(1-t)^{1/2}\, ,\, y_r(1-t)^{1/2})
     & \text{ when }y_r\le 0.
  \end{array} 
  \right.
\end{equation}
It is easily seen that $h_t$ preserves the equivalence relation, and is continuous.
The image of $h_1$ is 
the subset given by the points with $y_r\ge 0$, and $x_r=0$ whenever $y_r=0$. 
Since this subset is contractible in $\PP$ or $\PP_C$ (just `push inside' the points in the boundary), 
it follows that $\PP$ is contractible. Since $\PP_s$ and $\PP_C$ are invariant under $h_t$,
they are also contractible.
Notice that when $y_r=0$, the first coordinate of $h_t(x_r,y_r)$ is the opposite of that of $h_t(-x_r,y_r)$, 
so $h_t$ is well defined and continuous on $\mathbb{M}$, which is
therefore also contractible. Let us now look at $2\PP_C$.
First, we define another homotopy in $\PP_C$ 
that sends each collar to the line with $x_r$-coordinate $0$
by letting $\wh{h}_t(x_r,y_r)=h_t(x_r,y_r)$ if $y_r>0$ and
$$
  \wh{h}_t(x_r,y_r)= \bigl( x_r(1-t)^{1/2}, y_r\bigr) \text{ when }y_r<0.
$$ 
Let $\Gamma$ be the graph (and hence the CW-complex) given by two points $\{a,b\}$ linked by continuum many edges $e_r$
($r\in\R$). 
Map $\Gamma$ in $2\PP_C$ by $f$ as suggested on the picture below.

\begin{center}
\epsfig{figure=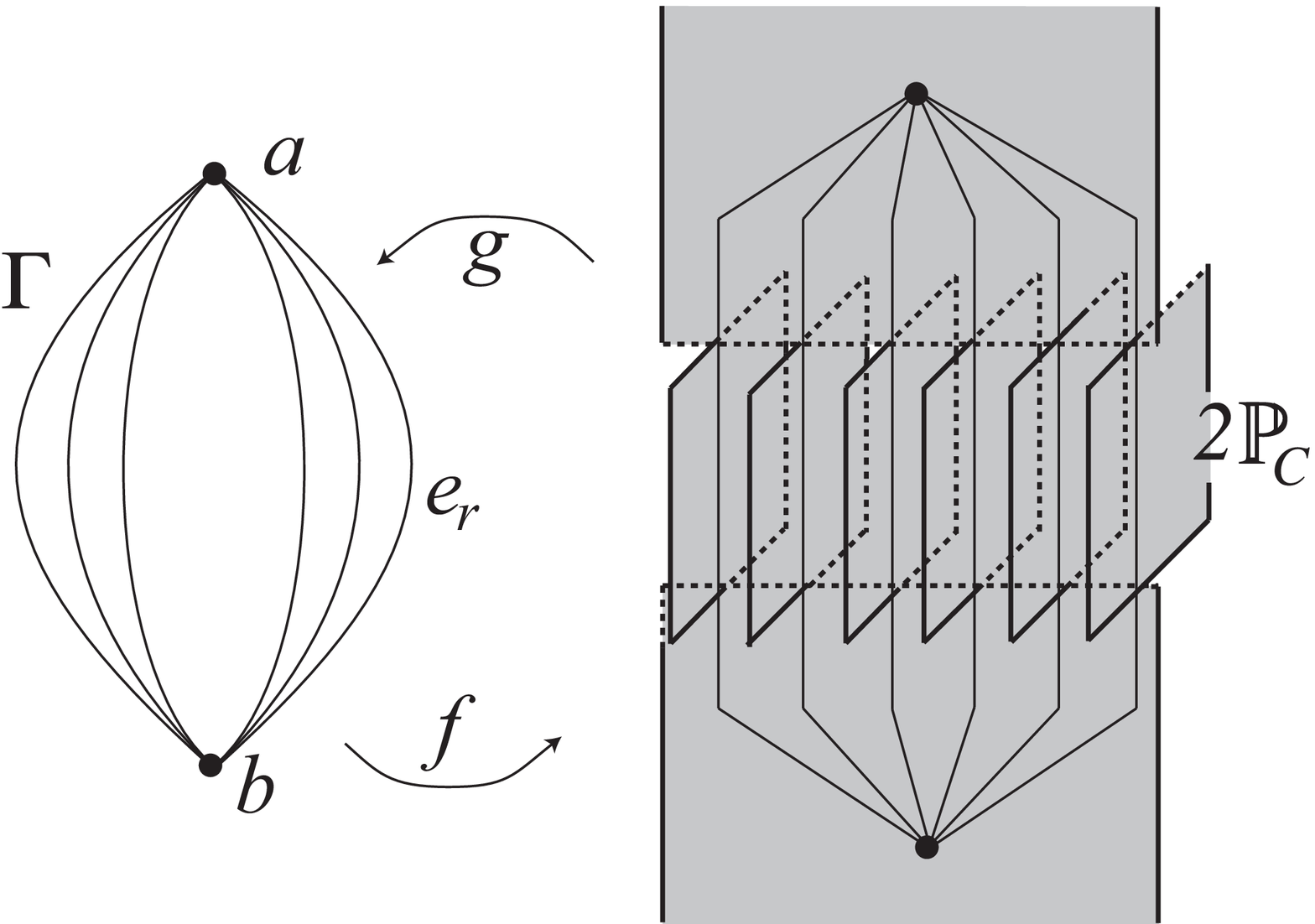, height=5cm}
\end{center}

It is not difficult (see \cite{GabardWouuuh_arxiv}) to show that $f$ is
a weak homotopy equivalence\footnote{i.e. induces isomorphisms on the $\pi_i$ for all $i\ge 0$} (though we do not need it for $2\PP_C$). 
In fact, defining $g:2\PP_C\to\Gamma$ as the map that sends the top and bottom regions respectively to $a$ and $b$ and
each `bridge' to the corresponding edge by projection on the $y_r$-factor, it is easy to check that
$g\circ f$ is homotopic to the identity in $\Gamma$.
But postcomposing $\wh{h}_t$ with a suitable homotopy in the top and bottom regions, we see that
$f\circ g$ is homotopic to the identity in ${2\PP_C}$ as well, so $2\PP_C$ is homotopy equivalent to $\Gamma$.

In the case of $2\PP_s$, a simple computation shows that $\pi_1(2\PP_s)=*_\R \Z$, the free group on a continuum of generators
(details can be found in \cite{GabardWouuuh}, though this is a direct application of van Kampen Theorem).
But $2\PP_s$ is separable, and
a separable space having the homotopy type of a CW-complex has the homotopy type of a countable CW-complex
(this is a direct consequence of Proposition A.1 p. 520 in \cite{Hatcher}),
and has thus
a countable fundamental group
(see \cite[Exercise 22, p. 360]{Hatcher},
it follows from the theorems on cellular approximations). 
This yields that $2\PP_s$ is is not homotopy equivalent to a CW-complex,
though a map $f:\Gamma\to 2\PP_s$ can be similarly defined and is again a weak homotopy equivalence.
\endproof

Now, the definitions.
If $M$ is a bordered manifold, write $\partial M$ for its boundary and $\text{\rm int}(M)$ for
its interior $M-\partial M$. To avoid confusion, we will always
specify when {\em topological} closures and interiors are taken.

Given an $n$-manifold $M$ and two bordered submanifolds $U,V\subset M$ whose interiors
are disjoint, we say that $B\subset S$ is a {\em bridge linking $U$ and $V$} if there is a homeomorphism
$\varphi:[a,b]\times N\to B$, where $N$ is an $n-1$-manifold, $\varphi(\{a\}\times N)$ is a boundary component of $U$, 
$\varphi(\{b\}\times N)$ is a boundary component of $V$, and $\varphi((a,b)\times N)$ is disjoint from $U\cup V$.
(If $a=b$, $[a,a]$ is understood as the singleton $\{a\}$.)
In words: we link together a boundary component of $U$ and one of $V$ inside $M$ by a collar of height $b-a$. If $a=b$, the bridge
is called infinitesimal.

Recall that $\mathbb{P}_s$ has boundary components $I_r$ which are real lines, indexed by $r\in\R$.
In $2\PP_s=\PP_s^0\cup\PP_s^1$, the two copies of $I_r$ are identified and thus form
an infinitesimal bridge linking $\PP_s^0$ and $\PP_s^1$.
In $2\mathbb{P}_C = \mathbb{P}_C^0 \cup \mathbb{P}_C^1$, we may take 
$U=\mathbb{P}_s^0\subset\mathbb{P}_C^0$, $V=\mathbb{P}_s^1\subset\mathbb{P}_C^1$, then the boundary components $I_r^0\subset \mathbb{P}_s^0$
is linked to $I_r^1\subset \mathbb{P}_s^1$ by a bridge which is the union of the two copies of $\{(x_r, y_r)\in\PP\,:\,-1\le y_r\le 0\}$, and are
thus {\em not} infinitesimal.

\begin{defi}
  Let $M$ be an $n$-dimensional manifold. A manifold-system of infinitesimal bridges (MSIB) in 
  $M$ consists of two 
  $n$-dimensional connected
  bordered submanifolds $U,V$
  such that $\text{\rm int}(U)\cap\text{\rm int}(V)=\varnothing$, $\text{\rm int}(U)\cup\text{\rm int}(V)$ 
  is metrizable, but 
  $\text{\rm int}(U)\cup\text{\rm int}(V)\cup (\partial U\cap\partial V)$ is not.
\end{defi}
 
In short, such a MSIB mimicks the behaviour of $2\PP_s$ (in the $2$-dimensional case of course).
Since  $\partial U\cap\partial V$ cannot be empty, 
each of its components is then an infinitesimal bridge.
The fullness condition below
says that we did not add `new bridges' in $U,V$ that are not intrinsic in $M$, by choosing badly $U$ and $V$ (see also Section \ref{sec3}).
For instance, a bad choice for $U,V$ in $2\PP_s$ would be to 
set $U=\PP_s^0 - D$, $V=\PP_s^1 -D$ where
$D$ is a set consisting of one point in each $I_r$. Thus, each $I_r$ gets cut artificially in two, adding new `non-intrinsic' bridges.

\begin{defi}
  An  MSIB $(U,V)$ in a manifold $M$ is full if $\partial U=\partial V$, and whenever $A$ is a boundary component of $\partial U$, 
  then $A$ is closed in $M$.
\end{defi}

\section{Proof of Theorem \ref{thmMSIB}}

The idea of the proof is the following. In the first step, we show that we can assume that an MSIB looks really like 
$2\PP_s$ 
(for instance, the boundary components are real lines except at most countably many circles). In the 
second step, we 
show that there is an uncountable family of loops in $U\cup V$ which are pairwise non-homotopic in the whole surface.
The theorem will follow using the separability of $\text{\rm int}(U)\cup \text{\rm int}(V)$.

\subsection{First step}

First, we show that we can assume that the surface does not contain $\omega_1$ (in the sense that there is a subset
homeomorphic to $\omega_1$). 
Recall that $\omega_1$ is the ordered set of all countable ordinals, and is a topological space 
with the order topology. It is probably the simplest example of a sequentially compact non-compact space.
It might seem that a bordered manifold with a metrizable interior can anyway not
contain $\omega_1$, but Nyikos found such a surface (whose boundary is the long ray $\LL_+$), see
\cite{Nyikos:1990}.

\begin{lemma}\label{CWnoomega_1}
  A manifold containing $\omega_1$ is not homotopy equivalent to a CW-complex.
\end{lemma}
\proof
  We shall abuse notation and write $\omega_1$ for the subset of $M$ homeomorphic to it.
  Let $K$ be a CW-complex and let $f:M\to K$ be an 
  homotopy equivalence with inverse $g:K\to M$.
  The image of $\omega_1\subset M$ under $f$ is sequentially compact in $K$ and thus compact (see the proof
  of Proposition A.1 p. 520 in \cite{Hatcher}). Then $g\circ f(\omega_1)\subset U$ where $U$ is a finite union of
  Euclidean open subsets of $M$. \\
  It is well known (and easy to prove) that any continuous map from $\omega_1$ to $\R^n$ is eventually constant, meaning 
  that there is some $\alpha\in\omega_1$ such that the map is constant on $[\alpha,\omega_1)\subset\omega_1$
  (see for instance \cite[Lemma 3.4 iii)]{Nyikos:1984}).
  This holds too for maps from $\omega_1$ to a finite union of Euclidean open subsets, for instance by embedding them is some $\R^n$. 
  (In fact, any continuous map $\omega_1\to X$ with $X$ first countable, Hausdorff and Lindel\"of is eventually constant, see Lemma 4.3 in
  \cite{BGGG}.)\\
  We now follow the proof of Proposition 7.1 in \cite{meszigueshom}.
  Suppose that
  there is a homotopy $h_t:M\to M$ with $h_0=g\circ f$ and $h_1=id$; $h_0\restrict{\omega_1}$ is thus eventually constant.
  Set 
  $$
    \tau = \inf\{t\,:\, h_s\restrict{\omega_1}\text{ is eventually constant }\forall s\le t\}.
  $$
  Then, $h_\tau$ is eventually constant, because either $\tau=0$, or there is some sequence $t_n\nearrow\tau$ ($n\in\omega$) with
  $h_{t_n}$ constant on $[\alpha_n,\omega_1)$ for some $\alpha_n$. Then, by continuity, there is some $x\in M$ such that
  $h_\tau([\beta,\omega_1))=\{x\}$ with $\beta=\sup_n \alpha_n$. 
  (Recall that a countable subset of $\omega_1$ is bounded, and thus possesses a supremum.)\\
  Since $h_1=id$, $\tau < 1$. 
  Let now $t_n\searrow \tau$, $n\in\omega$, 
  such that $h_{t_n}$ is not eventually constant, and let $V\ni x$ be an Euclidean neighborhood. 
  If $h_{t_n}$ sends $\omega_1$ inside $V$, then as above $h_{t_n}$ would be eventually constant, so there is some
  $\gamma_n\in\omega_1\subset M$, such that $h_{t_n}(\gamma_n)\in M-V$. We can choose $\gamma_n\ge \beta+1$.
  Taking a convergent subsequence for the $\gamma_n$ converging to some $\gamma > \beta$ yields
  that $h_\tau(\gamma)\in M-V$, contradicting the fact that
  $h_\tau([\beta,\omega_1))=\{x\}\subset V$.
  \\
  Thus, such a homotopy does not exist.
\endproof

By the classification of $1$-manifolds 
(a complete list is $\R,\mathbb{S}^1,\LL_+,\LL$, the latter two contain $\omega_1$) we have immediately:
\begin{lemma}\label{propbridges}
  Let $S$ be a surface with boundary that does not contain $\omega_1$. 
  Then each connected component
  of $\partial S$ is metrizable.
\end{lemma}

Lemma \ref{propbridges} is {\em false} in dimension $3$: the manifold with boundary $\PP_s\times[0,1]$ has a metrizable
interior but its boundary contains $\PP_s\times\{0\}$, which is non-metrizable and connected.
We now prove Proposition \ref{nocircles}, which says that a 
surface with boundary and metrizable interior has at most countably many circle boundaries. This proposition is not
really needed for the proof of Theorem \ref{thmMSIB}, but is quite simple and interesting in itself.
 
\proof[Proof of Proposition \ref{nocircles}]
  We may assume that $S$ has infinitely many boundary components that are circle.
  If $C$ is a circle boundary of $S$, 
  there is a neighborhood $N_C$ of $C$ in $S$ that intersects no other boundary component of $S$.
  By shrinking it if necessary, $N_C$ can be made homeomorphic to a cylinder $\mathbb{S}^1\times[0,1)$,
  and it retracts on $C$ (in $S$)
  (if needed, a proof can be found in \cite[Collaring Theorem 1.7.3, p. 35]{Rushing:Book}).
  There is thus an embedding $B_C$ of the circle in $N_C -  C\subset\text{\rm int}(S)$ 
  that `turns once' around $C$, obtained by pushing $C$ inside the collar. 
  We show that if $C\not= C'$, then $B_C$ and $B_{C'}$ cannot be homotopic (in $\text{\rm int}(S)$).
  Indeed, if there is a homotopy between them, then there is one between $C$ and $C'$ in 
  $\text{\rm int}(S)\cup C\cup C'$.
  Cover the image of the homotopy by finitely many Euclidean disks whose union is a 
  (metrizable) surface $S^-$
  with a boundary that 
  contains $C\cup C'$. 
  We sew a disk $D$ to $C'$, then $C$ becomes contractible in $S^-\cup D$, and thus bounds
  a $2$-disk in the same space by \cite[Theorem 1.7, p.85]{Epstein:1966} (for instance). 
  Removing $D$ again, we see that $S^-$ contains a cylinder with 
  boundaries $C,C'$, which is impossible since $S$ is connected.
  Thus $C$ and $C'$ cannot be homotopic in $S$.  The fundamental group of a metrizable manifold is at most countable
  (since it has the homotopy type of a countable CW-complex), so
  $\pi_1(\text{\rm int}(S))$ is at most countable, and there can be at most countably many circle boundary components.
\endproof

\begin{cor}\label{corlines}
  Let $U,V$ be a MSIB for some surface $S$ that does not contain $\omega_1$. 
  Then $\partial U\cap \partial V$ contains
  uncountably many connected
  components which are homeomorphic to intervals (open, closed or semi-open). If the MSIB is full, these intervals are 
  all open.
\end{cor}
\proof  
  By Proposition \ref{nocircles}, 
  $U$ and $V$ cannot have uncountably many boundary components that are circles.
  The intersection $A\cap B$ of connected components $A\subset\partial U$, $B\subset\partial V$ are all
  $1$-dimensional  
  manifolds (perhaps with boundary). Longlines and longrays are excluded since
  $\partial U\cap \partial V$ does not contain
  a copy of $\omega_1$. 
  Thus, $\partial U\cap\partial V$ 
  possesses uncountably many connected components. 
  Indeed, a manifold is metrizable if and only if it is Lindel\"of (recall that we assume connectedness), and if there
  are only countably many components in $\partial U\cap\partial V$, then 
  $U\cup V = \text{\rm int}(U)\cup\text{\rm int}(V)\cup(\partial U\cap\partial V)$ is Lindel\"of.\\
  The claim about full MSIB is immediate.
\endproof

\begin{lemma}\label{full-components}
  If $(U,V)$ is a full MSIB and $A$ is a component of $\partial U\cap\partial V$, then
  $A$ is a component of $S - (\text{\rm int}(U)\cup \text{\rm int}(V))$. 
  Moreover, any union of components of $\partial U\cap\partial V$ is closed in $S$.
\end{lemma}

\proof  
  Let $x\in A$, let 
  $U\subset A\cup \text{\rm int}(U)\cup \text{\rm int}(V)$
  be a small disk around $x$. Then its trace on $S - (\text{\rm int}(U)\cup \text{\rm int}(V))$ is
   included in $A$, which is therefore
  open. The `moreover' part follows from the fact 
  (see \cite[Collaring Theorem 1.7.3, p. 35]{Rushing:Book} again)
  that for each boundary component $A$  
  there is a neighborhood of $A$ (in $U\cup V$ for instance) whose closure does not     
  intersect any other component.
\endproof

\subsection{Second step}

Recall that a path in some space $X$ is a continuous map $p:[0,1]\to X$, and a loop is a 
path with $p(0)=p(1)$.
Inverses $p^{-1}$ and compositions $p_0p_1$ of paths and loops are defined as usual.
Our key lemma is the following:

\begin{lemma}\label{anewlemma}
   Let $M$ be a manifold and
   $\ell_\alpha:[0,1]\to M$ for $\alpha\in\kappa\ge\omega_1$ be loops that are pairwise non-homotopic in $M$.
   If $\cup_{\alpha\in\kappa}\ell_\alpha([0,1])\subset M$ is separable,
   $M$ does not have the homotopy type of a CW-complex. 
\end{lemma}

\proof
   Write $E=\cup_{\alpha\in\omega_1}\ell_\alpha([0,1])$.
   Suppose there is a homotopy equivalence $g:M\to K$, with $K$ a CW-complex.
   The image of a separable space is separable, thus $g(E)$ is separable.
   But a separable subset of a CW-complex is contained in a countable subcomplex,
  a fact already used in Proposition \ref{propop}.
  Call $K'\supset g(E)$ such a countable subcomplex.
  Since $g$ is a homotopy equivalence, the loops $g\circ\ell_\alpha$ in $K'\subset K$  are pairwise non-homotopic in $K$, and thus
  they are
  `even less' homotopic in $K'$.
  But a countable CW-complex has an at most countable fundamental group, a contradiction.
\endproof

For proving Theorem \ref{thmMSIB}, we shall exhibit a collection of loops that fulfill the assumptions of 
Lemma \ref{anewlemma}.
So, given a surface $S$ containing an MSIB $(U,V)$ and not containing $\omega_1$, 
let $\kappa$ be the
cardinal number of the connected components of $\partial U = \partial V$, 
and enumerate them as $A_\alpha$, $\alpha\in\kappa$ 
($\kappa\ge\omega_1$ by Corollary \ref{corlines}). 
By Proposition \ref{nocircles}, we can assume that there are at most countably many components of $\partial U$ that are circles,
that we may throw away,
and by Corollary \ref{corlines} and Lemma \ref{full-components} we obtain a full MSIB $U,V$ whose boundary components are lines. Call 
{\em flat} such a MSIB. \\
Fix two points $u,v$ in $\text{\rm int}(U)$ and $\text{\rm int}(V)$ respectively.
For each $\alpha\in\kappa$, let $p_\alpha$ be a path from $u$ to $v$ in
$A_\alpha\cup \text{\rm int}(U)\cup \text{\rm int}(V)$ (such a path exists since $A_\alpha\cup \text{\rm int}(U)\cup \text{\rm int}(V)$ is a subsurface).
We can choose $p_\alpha$ such that it intersects $A_\alpha$ in just one point.
Set $\ell_\alpha$ to be the loop $p_0 p_\alpha^{-1}$.
Since $\text{\rm int}(U)\cup \text{\rm int}(V)$ is a metrizable manifold, it is hereditarily separable (that is: any subset is separable), and 
since any point in the intersection of $\partial U$ and the image of $\ell_\alpha$ is the limit of a sequence of points
in $\text{\rm int}(U)\cup \text{\rm int}(V)$,
$\cup_{\alpha\in\kappa}\ell_\alpha([0,1])$ is separable.\\
If we can show that the $\ell_\alpha$'s are pairwise non-homotopic, $S$ satisfies
the assumptions of Lemma \ref{anewlemma}, and we are over.

\begin{lemma}\label{pi_1uctbl}
  Let $S$ be a surface with full and flat MSIB $(U,V)$.
  Let $\ell_\alpha=p_0p_\alpha^{-1}$ $\alpha\in\kappa\ge\omega_1$ be as above. 
  Then, the $\ell_\alpha$'s are pairwise non-homotopic.
\end{lemma}

The proof will follow from a reduction to the metrizable case, and the next lemma. 

\begin{lemma}\label{lemma_metrizable_surfaces}
  Let $T$ be a metrizable surface, $X,Y\subset T$ be connected subsurfaces with common boundary such that $\text{\rm int}(X)\cap \text{\rm int}(Y) =\varnothing$,
  and $\partial X = \partial Y$ is a finite disjoint union $\cup_{i=0,\dots,n}I_i$, with $n\ge 1$, $I_i$ closed in $T$ and homeomorphic to $\R$.
  Let $x\in \text{\rm int}(X), y\in \text{\rm int}(Y)$, and for $i=0,\dots,n$
  let $p_i$ be a path in $\text{\rm int}(X)\cup \text{\rm int}(Y)\cup I_i$ with $p(0)=a$, $p(1)=b$, such that 
  $p_i([0,1])\cap I_i$ is a single point.
  Set $\gamma_i = p_0 p_i^{-1}$, and let $\gamma=\lambda_0\cdots\lambda_k$, with $\lambda_j=\gamma_i$ or $\lambda_j=\gamma_i^{-1}$ 
  for $i\in\{1,\dots,n\}$ be a reduced
  loop. (By this we mean that corresponding word is reduced.) \\
  Then, $\gamma$ is not nullhomotopic as a loop in $T$.
\end{lemma}

\proof[Proof of Lemma \ref{pi_1uctbl}]
  In fact, we show a little bit more:
  any nonempty reduced word with letters in $\{\ell_\alpha,\ell_\alpha^{-1}\,:\,1\le\alpha<\kappa\}$
  is not nullhomotopic as a loop in $S$.
  Suppose thus that a nontrivial product $\ell$ of these $\ell_{\alpha}$'s is nullhomotopic in $S$.
  Call $H$ the image of the homotopy.
  By compactness, there is a finite number 
  $A_{\alpha_{n+1}},\dots, A_{\alpha_m}$ of 
  components of $\partial U=\partial V$ 
  different from $A_{\alpha_{0}},\dots, A_{\alpha_n}$ that intersect $H$. 
  Cover $H$ by Euclidean open sets, by compactness there is a
  finite subcover, and 
  let $T$ be the union of the members of this subcover, $X=U\cap T_0$ and $Y=V\cap T$.
  Then, $A_{\alpha_i}\cap T$ is closed in $T$ and homeomorphic to some disjoint union of open intervals. By removing all the components
  that do not intersect $H$, we can assume that this union is finite.
  Letting $I_i$ be the component of $A_{\alpha_i}\cap T$ that intersects $p_{\alpha_i}$, the assumptions of Lemma 
  \ref{lemma_metrizable_surfaces} are fulfilled, so $\ell$ is not homotopic to a constant map, a contradiction.
\endproof

To finish, we just need to prove Lemma \ref{lemma_metrizable_surfaces}, a task we tackle now.

\proof[Proof of Lemma \ref{lemma_metrizable_surfaces}]
Let $\pi:\wt{T}\to T$ be the universal covering of $T$. Since $T$ is not compact 
(it contains closed non-compact subsets, namely the $I_i$'s), 
so is $\wt{T}$, which is thus homeomorphic to $\R^2$ (see \cite[Corollary 1.8]{Epstein:1966}).
If $E$ is a subset of $T$, we denote by $\wt{E}$, $\wt{E}'$, $\wt{E}''$ different lifts of $E$ in $\wt{T}$.
(We slightly abuse language here as we should consider lifts of the inclusion map.)
Since $I_i$ is contractible, each $\wt{I}_i$  is a 
copy of $I_i\approx\R$ and is closed in $\wt{T}$. \\
Sch\"onfliess theorem is equivalent to the assertion that any closed copy of $\R$ in $\R^2$ disconnects the space in two components both homeomorphic
to $\R^2$.
(To see this, take the one point compactification of $\R^2$ which is $\mathbb{S}^2$. 
The copy of $\R$ gets closed by the added point, and thus is homeomorphic to
$\mathbb{S}^1$. Then apply Sch\"onfliess. For the converse, take out a point of the embedded circle.)
So, $\wt{I_i}$ disconnects $\wt{T}$. But in fact, since $\wt{I_i}\cap\wt{I_j}=\varnothing$ whenever $i\not= j$
and $\wt{I_i}\cap\wt{I_i}'=\varnothing$, given a connected component $C$ of 
$\wt{T} - (\wt{I_{i_1}}\cup\cdots\cup\wt{I_{i_k}})$, if $\wt{I_i}$ intersects $C$ then it disconnects it. We refer to those 
$\wt{I_i}$'s as {\em $i$-fences}.\\
Fix a lift $\wt{x}$ of $x$, and let $\wt{\gamma}$ be the lift of $\gamma$ with $\wt{\gamma}(0)=\wt{x}$.
By the definition of the path $p_i$, $\wt{\gamma}$ `goes through' the fences in the order prescribed by the word $\gamma$.
For instance, if $\gamma=\gamma_i=p_0p_i^{-1}$, then $\wt{\gamma}$ goes through a $0$-fence and then an $i$-fence. Moreover, $\gamma$ intersects
only those two among all the fences.
\\
Set $m(\gamma)$ to be the smallest number of $i$-fences with $i\ge 1$ that a path linking $\wt{\gamma}(0)=\wt{x}$ to $\wt{\gamma}(1)$ has to cross.

\begin{lemma}\label{sublemma_surfaces}
   Let $\gamma = \lambda_{i_0}\cdots \lambda_{i_k}$ be as above (in particular, it is reduced as a word in the $\gamma_i,\gamma_i^{-1}$, $i=1,\dots,n$). 
   Then $m(\gamma)=k+1$.
\end{lemma}

\proof 
   By induction.
   If $\gamma=\gamma_i$ or $\gamma_i^{-1}$, then as said before
   it crosses two fences which are disjoint and disconnect $\wt{T}$.
   Thus, a path joining $\wt{\gamma_i}(1)$ to $\wt{x}$ must cross at least these two fences, 
   one of which is an $i$-fence with $i\ge 1$,
   so $m(\gamma)=m(\gamma_i)=1$.\\
   Suppose now that the lemma holds for $\gamma=\lambda_{i_{0}}\cdots\lambda_{i_{k-1}}$, and let
   $\gamma'=\gamma\cdot\lambda_{i_{k}}$.
   We have 
   $\wt{\gamma}(1)=\wt{x}'$, and $\wt{\gamma}'(1)=\wt{x}''$. Suppose that $\lambda_{i_{k}}=\gamma_i$, $i\ge 1$.
   If $\lambda_{i_{k-1}}=\gamma_{i'}$, then $\wt{\gamma}$ crosses first a $0$-fence and then an $i'$-fence.
   If $\lambda_{i_{k-1}}=\gamma_{i'}^{-1}$, it crosses first an $i'$-fence and then a $0$-fence. 
   Since $i'\not= i$ in the latter case (the word is reduced), we conclude in both cases that
   a path linking $\wt{x}'$ to $\wt{x}''$ crosses at least an 
   $i'$-fence\footnote{Actually, if $\gamma_{i_{k-1}}=\gamma_{i_k}^{-1}$, then $m(\gamma')=m(\gamma)-1$.}.
   Moreover, $\wt{\gamma}$ crosses exactly $k$ fences, which is the minimal number 
   for a path linking $\wt{x}$ to $\wt{x}'$ by induction.
   It follows by the above `disconnecting remarks' that $\wt{\gamma'}$ crosses the minimal number of $i$-fences with $i\ge 1$ among the paths 
   linking $\wt{x}$ to $\wt{x}''$, thus $m(\gamma')=k + 1$, as wished.
   The case $\lambda_{i_{k}}=\gamma_i^{-1}$, $i\ge 1$, is treated exactly the same.
\endproof

  This enables to finish the proof of Lemma \ref{lemma_metrizable_surfaces}: if the loop $\gamma$ is contractible, then
  $\wt{\gamma}(0)=\wt{\gamma}(1)$, which is impossible since they are separated by $k+1$ fences.
\endproof

\begin{rem}
  In fact, the proofs of Lemmas \ref{pi_1uctbl}--\ref{lemma_metrizable_surfaces} work as well if some of the $A_\alpha$ 
  and $I_i$ are circles. Indeed, any circle in
  $\partial U=\partial V$ cannot bound a $2$-disk by connectedness of $U$ and $V$, and thus its lifts have 
  the same `disconnection properties'. 
  Proposition \ref{nocircles} is therefore not needed in the proof of Theorem \ref{thmMSIB}.
\end{rem}

\subsection{Other related results}

In this subsection we give some results that, though not needed for the proof of Theorem \ref{thmMSIB}, are somewhat related. 
(Actually, we used them in a first attempt at the proof.)
First, we show that the loops $\ell_\alpha$ yield a weak homotopy equivalence between $U\cup V$ and a graph (if the MSIB is flat).

\begin{lemma}\label{whe1} 
  Let $S$ be a surface with a full flat MSIB $(U,V)$.
  Then
  there is a weak homotopy equivalence between a graph $\Gamma$ and $S_0$. 
\end{lemma}

\proof
  We write $S_0=U\cup V$, $S_0$ is a separable submanifold.
  ${\rm int}(U)$ and ${\rm int}(V)$, being open metrizable surfaces, have the homotopy type of an at most   
  countable 
  wedge of circles $W_U,W_V$ with vertex $a$ and $b$ respectively. Recall that 
  $A_\alpha$ ($\alpha\in\kappa$) are the connected components of $\partial U=\partial V$.
  Build $\Gamma$ adding edges $e_\alpha$, $\alpha\in\kappa$,
  linking $W_U$ to $W_V$.
  The map $f:\Gamma\to S_0$ defined by sending the vertices of $\Gamma$ to $u$ and $v$, the loops 
  in $W_U,W_V$ to the generators of the $\pi_1$ of $\text{\rm int}(U),\text{\rm int}(V)$ respectively and each edge
  $e_\alpha$ to the path $p_\alpha$ induces an isomorphism of the fundamental groups (see \cite{GabardWouuuh_arxiv} for details, this is essentially 
  a direct application of van Kampen's Theorem).\\
  It is a standard fact that a metrizable non-compact surface has vanishing $n$-homotopy groups for $n\ge 2$, since it 
  which is homotopy equivalent to an at most countable
  wedge of circles. 
  Given $f:\mathbb{S}^n\to S$, we can cover its image by finitely many Euclidean charts whose union forms a non-compact metrizable surface.
  Thus, $\pi_n(S)=0$ for $n\ge 2$, and
  $f$ is therefore a weak homotopy equivalence.
\endproof

In fact we can obtain a little more:

\begin{lemma}\label{whe2} 
  Let $S$ be a surface with flat MSIB $(U,V)$, $S_0=U\cup V$.
  There is a CW-approximation $(M,\Gamma)$ of the pair $(S,S_0)$ such that $\Gamma$ is the above mentioned graph. 
\end{lemma}

\proof
  As in \cite[ex. 4.15 p.353]{Hatcher}, we first CW-approximate $S_0$ by $f_0:\Gamma\to S_0$ using Lemma   \ref{whe1}, and then 
  build a CW-approximation of the pair $(N,\Gamma)$ where $N$ is the mapping cylinder of the composition of
  the inclusion $S_0\to S$ and $f_0$. 
\endproof

\section{Possible generalisations}\label{sec3}
  
Can the definition of full MSIB be weakened such that Theorem \ref{thmMSIB} still holds~?
For instance, can we prove something with the following definition:
\begin{defi}
  Let $M$ be a manifold. A system of infinitesimal bridges (SIB) in $M$ consists of two disjoint connected open   
  sets $U,V$
  such that $U\cup V$ is metrizable but $U\cup V \cup (\wb{U}\cap\wb{V})$ is not.
\end{defi}

Of course, since topological boundaries can be quite nasty, 
it is a priori not clear whether the nice `disconnection properties' we used in Lemma \ref{lemma_metrizable_surfaces}
are still valid.

Can we relax somehow the definition of fullness~? 
for instance, does this weaker definition suffices~?

\begin{defi}
  An  MSIB $(U,V)$ in a manifold $M$ is weakly full if whenever $x,y$ belong to two different connected  
 components of 
  $\partial U\cap\partial V$, then they belong to two different connected components of   
  $\partial U$ and
  of $\partial V$.
\end{defi}

The idea here is that you do not add artificial bridges because 
the boundary components of $U$ and $V$ do not intersect well. For instance, 
consider as before $2\PP_s$, set $U$ to be one of the copies of $\PP_s$ in $2\PP_s$, and $V$ to be the other
to which a point in each boundary component has been removed. You have thus `artificially'
cut in two each bridge, and there are generators of $\pi_1(U\cup V)$ that vanish
in $\pi_1(2\PP_s)$. 
This weak definition of fullness prevents this kind of bad
intersections. But it is not strong enough to ensure that 
Lemma \ref{pi_1uctbl} works.
Take for instance the surface $S$ of Nyikos alluded above: 
$S$ has metrizable
interior and boundary $\LL_+$, and consider the double $2S$, it 
has trivial homotopy groups.
But if we denote the two copies of $\text{\rm int}(S)$ by $S_1,S_2$ and let $U=S_1\cup(\LL_+ - \omega_1)$,
$V=S_2\cup(\LL_+ -\omega_1)$, then $(U,V)$ form a MSIB that is weakly full. 
(Of course, $2S$ does not contradict Theorem \ref{thmMSIB}, because of Lemma \ref{CWnoomega_1}.)

More generally, there is a variety of 
questions as whether a surface containing a [weakly full] 
(M)SIB always contains a [full] (M)SIB, taking or not the ()'s and []'s parts giving the various possibilities. 

\begin{prob}
  What about higher dimensions ?
\end{prob}

The given proof would have some gaps in higher dimension, for instance
  Lemma \ref{propbridges} is false
in dimension $\ge 3$, 
and we used properties of the universal covering that are special for dimension $2$ in the proof
of Lemma \ref{lemma_metrizable_surfaces}.
Related is the following problem:

\begin{prob}
  Given a manifold $M$ containing a submanifold that is not homotopy equivalent to a CW-complex,
  can we conclude that $M$ is not homotopy equivalent to a CW-complex~?
\end{prob}

To finish, another general problem in dimension $2$:

\begin{prob}
  Let $S$  be a non-compact surface (metrizable or not).
  Is $\pi_1(S)$ a free group~?
\end{prob}

This is of course well known for metrizable surfaces. 

{\em Acknowledgements.}
For useful conversations, I would like to thank David Gauld, Sina Greenwood, and
David Cimasoni. I thank Alexandre Gabard for his numerous suggestions.


\begin{thebibliography}{10}

\bibitem{meszigueshom}
M.~Baillif.
\newblock The homotopy classes of continuous maps between some non-metrizable
  manifolds.
\newblock {\em Top. App.}, 148(1--3):39--53, 2005.

\bibitem{BGGG}
M.~Baillif, A.~Gabard, and D.~Gauld.
\newblock Foliations on non metrizable manifolds: absorbtion by a {C}antor
  black hole.
\newblock submitted, 2008.

\bibitem{meszigues+Nyikos}
M.~Baillif and P.~Nyikos.
\newblock Various homotopy types in type {I} spaces and manifolds.
\newblock in preparation, 2010.

\bibitem{CalabiRosenlicht}
E.~Calabi and M.~Rosenlicht.
\newblock Complex analytic manifolds without countable base.
\newblock {\em Proc. Amer. Math. Soc.}, 4:335--340, 1953.

\bibitem{Epstein:1966}
D.~B.~A. Epstein.
\newblock Curves on 2-manifolds and isotopies.
\newblock {\em Acta Math.}, 115:83--107, 1966.

\bibitem{GabardWouuuh_arxiv}
A.~Gabard.
\newblock A separable manifold failing to have the homotopy type of a
  {CW}-complex (first version).
\newblock eprint arXiv:math/0609665, 2006.

\bibitem{GabardWouuuh}
A.~Gabard.
\newblock A separable manifold failing to have the homotopy type of a
  {CW}-complex.
\newblock {\em Arch. Math.}, 90:267--274, 2008.

\bibitem{Hatcher}
A.~Hatcher.
\newblock {\em Algebraic topology}.
\newblock Cambridge University Press, Cambridge, 2002.
\newblock (A free electronic version is available).

\bibitem{MilnorCW}
J.~Milnor.
\newblock On spaces having the homotopy of a {CW}-complex.
\newblock {\em Trans. Am. Math. Soc.}, 90:272--280, 1958.

\bibitem{Nyikos:1984}
P.~Nyikos.
\newblock The theory of nonmetrizable manifolds.
\newblock In {\em Handbook of {S}et-{T}heoretic {T}opology}, pages 633--684.
  North-Holland, Amsterdam, 1984.

\bibitem{Nyikos:1990}
P.~Nyikos.
\newblock {\em On first countable, countably compact spaces III: The problem of
  obtaining separable non- compact examples,}, pages 128Ð--161.
\newblock North-Holland, Amsterdam, Amsterdam, 1990.

\bibitem{Rushing:Book}
T.~B. Rushing.
\newblock {\em Topological Embeddings}.
\newblock Number~52 in Pure and Applied Math. Academic Press, New York and
  London, 1973.

\end{thebibliography}

\end{document}